\documentclass[twoside]{article}

\usepackage{PRIMEarxiv}

\usepackage[utf8]{inputenc} 
\usepackage[T1]{fontenc}    
\usepackage{hyperref}       
\usepackage{float}
\usepackage{url}            
\usepackage{booktabs}       
\usepackage{amsfonts}       
\usepackage{microtype}      
\usepackage{lipsum}
\usepackage{fancyhdr}       
\usepackage{graphicx}       
\graphicspath{{media/}}     
\usepackage{tabularx}
\usepackage{amssymb}

\usepackage{amsmath}

\DeclareMathOperator*{\argmin}{arg\,min}

\title{Volterra black-box models identification methods: direct collocation vs least squares
}

\author{
  Denis Sidorov\\
  Applied Mathematics Department \\ 
  Melentiev Energy Systems Institute SB RAS \\
  Irkutsk, Russia\\
  \texttt{dsidorov@isem.irk.ru} \\
   \And
  Aleksandr Tynda \\
  Higher and Applied Mathematics Department \\
  Penza State University \\
  Penza, Russia\\
  \texttt{tyndaan@mail.ru} \\
   \And
  Vladislav Muratov \\
  Institute of Mathematics and Information Technologies \\
  Irkutsk State University \\
  Irkutsk, Russia\\
  \texttt{muratov428@gmail.com} \\
  \And
  Eugeny Yanitsky \\
  Intermediate Radio Frequency Lab \\
  Huawei Russian Research Institute \\
  Moscow, Russia\\
  \texttt{yanitsky.eugeny@huawei.com} \\
}

\begin{document}
\maketitle

\begin{abstract}
The Volterra integral-functional series is the classic approach for nonlinear black box dynamical systems modeling. It is widely employed in many domains including radiophysics, aerodynamics, electronic and electrical engineering and many other. Identifying the time-varying functional parameters, also known as Volterra kernels, poses a difficulty due to the curse of dimensionality. This refers to the exponential growth in the number of model parameters as the complexity of the input-output response increases. The least squares method (LSM) is widely acknowledged as the standard approach for tackling the issue of identifying parameters. Unfortunately, the LSM suffers with many  drawbacks such as the sensitivity to outliers causing biased estimation, multicollinearity,  overfitting and inefficiency with large datasets.
This paper presents alternative approach based on direct estimation of the Volterra kernels using the collocation method. Two model examples are studied.
It is found that the collocation method presents a promising alternative for optimization, surpassing the traditional least squares method when it comes to the Volterra kernels identification including the case when input and output signals suffer from considerable measurement errors.
\end{abstract}


\section{Introduction}

\begin{figure}[H]
\centering
\includegraphics[width=0.41\linewidth]{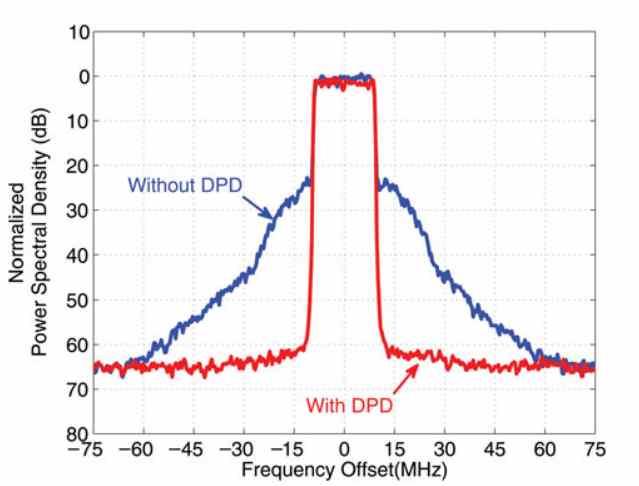}
\caption{Power spectrum density \cite{DPD_ZHU}} \label{FFig1}
\end{figure}

At the current stage of development of wireless technologies like 5G/6G communication system networks based on antenna arrays with digital beam forming (Massive Multiple Input Multiple Output system), it is impossible to do without such digital signal processing algorithms as digital correction of the nonlinear distortion DPD (Digital Predistortion). Nonlinear distortions of the signal occurring inside the transceiver path strongly distort the spectrum of this signal, as shown in Fig. \ref{FFig1}, where its shown in blue, and the main signal is red color respectively.
 
\begin{figure}[H]
\centering
\includegraphics[width=0.6\linewidth]{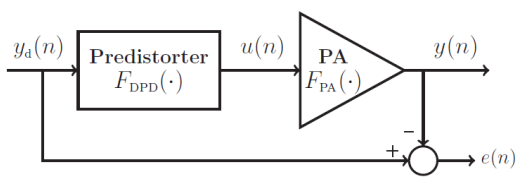}
\caption{Digital correction scheme} \label{Ffig2}
\end{figure}

However, the international wireless standards like 3GPP, ETSI 
impose strict requirements on the spectral power of the radiated signal. The use of digital nonlinear distortion correction algorithms allows to meet the requirements of standards and at the same time positively affect the overall efficiency, that is, the energy consumption of the entire signal receiving and transmitting system.
There are different approaches to the implementation of such algorithms, both purely digital and analog and mixed. One of them, a purely mathematical approach to the description of nonlinear distortions, we will describe below. However, Let's consider the general statement of the problem of digital correction (DPD) with the following structure of the some model of correction as shown in Fig. \ref{Ffig2}.
Here $F_{DPD}(.)$ is a nonlinear operator reflecting the essence of nonlinear correction and imagine it as some function dependent on parameters $\vec{W}=[w_{1},...,w_{p}]^T,\vec{W}\in \mathbf{C}^p$. $F_{PA}(.)$ is a nonlinear operator identified with a nonlinear device and generating some complex vector
$\vec{Y}=[y_{1},...,y_{n}]^T,\vec{Y}\in \mathbf{C}^n$ and also define some vector from a complex field of numbers $\vec{Y_d}=[y_{d,1},...,y_{d,n}]^T,\vec{Y_d}\in \mathbf{C}^n$
on which the operator $F_{DPD}(.)$  depends. Under the error $E\in \mathbf{C}^n$ we will understand the difference between vectors $Y$ and $Y_d$

$$\vec{E} = \vec{Y_d} - \vec{Y}.$$
Then we can formulate the requirements for the definition of parameters $\vec{W}$ as follows:
$\vec{\omega }=\argmin_W\left\| E \right\|^2$, where $\left\| . \right\|$ is Euclidean norm. Considering  $Y=F_{PA}(F_{DPD}(Y_d)),$  the above introduced expression can be rewritten as 

$$\vec{\omega }=\argmin_W\left\| Y_d - F_{PA}(F_{DPD}(Y_d)) \right\|^2.$$
This equation will be task of DPD (Digital Predistortion).
Here we can highlight several important sub-tasks, which in themselves are quite complex both theoretically and computationally:

a)  Since we have formulated, in fact, the problem of approximation of a function, we need to derive the analytical regression dependence $F_{DPD}(.)$ on the parameters $\vec{W}$. How this function is defined will depend on the quality of the correction of nonlinear distortions;

b) A procedure of the searching the parameters $\vec{W}$ is a classical optimization problem, which is a linear or nonlinear regression with respect to the parameters $\vec{W}$. Finding efficient methods of convex or non-convex optimization is one of the major problems in this problem;

c) Compression of a function $F_{DPD}(.)$, i.e. reducing its computational complexity.

One of the methods to solve the problem a) for DPD task is the Volterra functional series. And it is also the conventional tool to characterize 
the complex nonlinear dynamics in various fields including the radiophysics, mechanical engineering, electronic and electrical engineering, energy sciences (here readers may refer e.g. to review \cite{VSrev}).
Volterra series are widely employed to represent the input-output relationship of nonlinear dynamical systems with memory.
Volterra power series are among the best-understood nonlinear system representations in signal processing. Such integral functional series (also called Fr\'{e}chet-Volterra series) \eqref{VS} 
 \begin{equation}\label{VS}
  y(t)= F(x(t)):=\int\limits_{0}^{t}K_1(s)x(t-s)\,ds + \int\limits_{0}^{t}\int\limits_{0}^{t}K_2(s_1,s_2) x(t-s_1)x(t-s_2)\,ds_1ds_2 +
  \dots   
\end{equation}
$$ \dots +\int\limits_{0}^{t}\int\limits_{0}^{t} \dots \int\limits_{0}^{t} K_n(s_1,s_2, \dots, s_n)x(t-s_1)x(t-s_2)\dots x(t-s_n)  \, ds_1 ds_2 \dots ds_n + \dots \; t\in[0,T]$$
were first proposed by Maurice Fr\'{e}chet for a continuous nonlinear dynamical systems representation \cite{fre, volt}.
Here readers may also refer to overview \cite{vain} and monograph \cite{sid} for more details on relevant 
Lyapunov--Liechtenstein operator and  Lyapunov -- Schmidt methods in the theory of non-linear equations.

The role of a reproducing kernel Hilbert space in development of a unifying view of the Volterra theory and polynomial kernel regression is presented in \cite{fra}.

In (1) $x(t)$ is input signal and $y(t)$ is output  
of a single input single output (SISO) nonlinear system and $K_n(s_1,s_2, \dots, s_n)$ are the multidimensional Volterra kernels (or transfer functions) to be identified based on nonlinear system's response $y(t)$ as a reaction on input $x(t)$ (Fig.~3). It is to be noted, that for the basic case $n=1$
we have a conventional Finite Impulse Response (FIR) linear model which optimal in the least-squares sense. 
\begin{center}
\begin{picture}(300,100)(-40,0)
\put(50,80){\line(1,0){120}}
\put(20,70){\makebox(0,0){$x(t)$}}
\put(110,60){\makebox(0,0){
 $F(x(t))$ 
}}
\put(200,70){\makebox(0,0){$y(t)$}}
\put(0,60){\vector(1,0){45}}
\put(175,60){\vector(1,0){40}}
\put(50,80){\line(0,-1){40}}
\put(50,40){\line(1,0){120}}
\put(170,80){\line(0,-1){40}}
\put(110,15){\makebox(0,0){{\bf Figure 3}. Behavioral modeling of the black box system}}
\end{picture}
\end{center}

\setcounter{figure}{3}

Fr\'{e}chet theorem \cite{fre}  generalises the 
famous Weierstrass approximation theorem  which characterizes the set of continuous functions on a compact interval via uniform approximation by algebraic polynomials.

Power series (1) characterize the stationary dynamical systems. 
Stationarity here means that a transfer functions do not vary during the transient process as $t\in[0,T]$. More general power series \eqref{VS1}
models nonstationary dynamics when transfer functions
depend explicitly on time $t$
 \begin{equation}\label{VS1}
  y(t)= \int\limits_{0}^{t}K_1(t,s)x(s)\,ds + \int\limits_{0}^{t}\int\limits_{0}^{t}K_2(t,s_1,s_2)x(s_1)x(s_2)\,ds_1ds_2 +
  \dots   
\end{equation}
$$ \dots +\int\limits_{0}^{t}\int\limits_{0}^{t} \dots \int\limits_{0}^{t} K_n(t,s_1,s_2, \dots, s_n)x(s_1)x(s_2)\dots x(s_n)  \, ds_1 ds_2 \dots ds_n + \dots \; t\in[0,T].$$

The Volterra series is essential tool of mathematical 
modeling the nonlinear dynamical systems appearing in   
digital pre-distortion (DPD) iterative process \cite{dpd}. DPD as we described before is an important part of the digital signal processing algorithms used in transmitters and receivers. 
Several methods have been studied for DPD, with Volterra series-based methods being popular due to their ease of implementation and the straightforward interpretation of their nonlinear terms. The key issue with Volterra series is the curse of dimension: as the order of the series increases, the number of terms involved in the expansion grows exponentially, making it computationally demanding.
From other hand, estimating the functional coefficients (Volterra kernels) of the Volterra integral functional series can be challenging. It often considered 
in its discrete form and requires a significant amount of data and complex optimization algorithms to find the best fit for the model coefficients. Alternative approach based on problem reduction to multi-dimensional integral 
equations solution \cite{ref1, ref2} needs special probe signals design.

In present paper the alternative approach for the Volterra kernels 
identification is 
proposed using the direct collocation method. The results
are compared 
with the conventional least squares method (LSM) widely employed 
for Volterra series identification problem in telecommunication 
domain.

   The rest of the paper is structured as follows: The subsequent section provides the problem statement. Section 3 focuses on the collocation method. Section 4 carries out computational experiments with LSM, while section 5 discusses concluding remarks and future work.

\section{Identification problem statement}

Let us consider the following segment of the truncated Volterra series \eqref{VS} for \(n=2\)
\begin{equation}\label{VSeries-1}
  y(t)=\int\limits_{0}^{t}K_1(s)x(t-s)\,ds + \int\limits_{0}^{t}\int\limits_{0}^{t}K_2(s_1,s_2)x(t-s_1)x(t-s_2)\,ds_1 ds_2, \; t\in[0,T].
\end{equation}

Our current problem in this section is to determine the kernels \(K_1(s)\) and \(K_2(s_1,s_2)\)  by a known input and output pair \(\bigl(x(t),y(t)\bigr)\).

In contrast to the linear case \(n=1\), when it is sufficient to specify a single pair \(\bigl(x(t), y(t)\bigr)\) to determine the kernel \(K_1(s)\), in the nonlinear case \(n=2\),  for the unique identification of the two-dimensional kernel \(K_2(s_1, s_2)\), it is necessary to specify a two-dimensional continuum of equalities. This means that problem \eqref{VSeries-2} has an infinite set of solutions. 

 \newcounter{remark}
 \setcounter{remark}{0}
 \newtheorem{Remark}[remark]{Remark}
\renewcommand{\theremark}{A\arabic{remark}}

\begin{Remark} 
It should be noted that if we consider this problem as an integral equation with two unknown functions \(K_1(s)\) and \(K_2(s_1,s_2)\), then this problem is essentially ill-posed. 
There are an infinite number of solutions, this problem is insufficiently defined. In this regard, no classical numerical methods designed for integral equations are applicable in this case. 
And as a result, there are no any attempts to solve the problem in this form in the literature.
\end{Remark}

\begin{Remark}
A fundamentally different situation takes place in the problem of determining an unknown input signal \(x(t)\) with a known output signal \(y(t)\) after kernels identification.
It is to be noted that in this case we have the problem of nonlinear Volterra integral equations' solution.
Here readers may refer to
sec. 9 in book \cite{ref2}, papers \cite{ap1}, \cite{sidsid} and references therein regarding the Kantorovich main solutions and blow up phenomenon. 
\end{Remark}

Within the framework of this paper, from a practical point of view, we will be satisfied with any pair of approximately found kernels \(\widetilde{K}_1(s)\) and \(\widetilde{K}_2(s_1, s_2)\) that provides a sufficiently small residual norm
\begin{equation}\label{VSeries-2}
  \varepsilon=\max\limits_{t\in[0,T]}\left|y(t)-\int\limits_{0}^{t}\widetilde{K}_1(s)x(t-s)\,ds - \int\limits_{0}^{t}\int\limits_{0}^{t}\widetilde{K}_2(s_1,s_2)x(t-s_1)x(t-s_2)\,ds_1 ds_2\right|.
\end{equation}

Denote by \(B_i(t),\;i=0,1,2,\ldots,\) the basis functions forming a complete orthogonal system of functions on the segment \([0,T]\).

We look for an approximate solution of the problem \eqref{VSeries-1} in the form of segments of series of expansions according to the selected system of basis functions
\begin{equation}\label{VSeries-3}
  \widetilde{K}_{1,m}(s)=\sum\limits_{i=0}^{m-1}A_iB_i(s), \quad   \widetilde{K}_{2,m_1,m_2}(s_1,s_2)=\sum\limits_{i=0}^{m_1-1}\sum\limits_{j=0}^{m_2-1}C_{ij}B_i(s_1)B_j(s_2).
\end{equation}

\section{Collocation method}
Collocation-type methods are widely used in the discretization of various kinds of integro-functional equations \cite{Brunner-2004}. With sufficiently good accuracy and stability, they are also computationally less expensive in comparison with projection methods of the Galerkin type requiring additional integration \cite{Polyanin}.

In order to determine the unknown coefficients \(A_i\) and \(C_{ij}\), we introduce a uniform grid of nodes
\begin{equation}\label{VSeries-4}
   t_k\in[0,T], \; k=0,1,\ldots,N,
\end{equation}
where $N + 1$ is number of nodes.

Substitute \eqref{VSeries-3} in \eqref{VSeries-1} and then demand that the equalities be fulfilled at the points \eqref{VSeries-4}
\begin{equation}\label{VSeries-5}
  y(t_k)=\int\limits_{0}^{t_k}\widetilde{K}_{1,m}(s)x(t_k-s)\,ds + \int\limits_{0}^{t_k}\int\limits_{0}^{t_k}\widetilde{K}_{2,m_1,m_2}(s_1,s_2)x(t_k-s_1)x(t_k-s_2)\,ds_1ds_2, \; k=\overline{0,N}.
\end{equation}

Denote for a simplicity \(y(t_k)=y_k\), and transform the last equalities as follows
\begin{equation}\label{VSeries-6}
  y_k=\sum\limits_{i=0}^{m-1}A_i\int\limits_{0}^{t_k}B_i(s)x(t_k-s)\,ds + \sum\limits_{i=0}^{m_1-1}\sum\limits_{j=0}^{m_2-1}C_{ij}\int\limits_{0}^{t_k}\int\limits_{0}^{t_k}B_i(s_1)B_j(s_2)x(t_k-s_1)x(t_k-s_2)\,ds_1ds_2.
\end{equation}

As a system of basis functions \(B_i(t), \; i=0,1,\ldots\), we choose Chebyshev polynomials of the first kind
\begin{equation}\label{VSeries-9}
   T_0(t)=1,\; T_1(t)=t, \; T_{i+1}(t)=2tT_i(t)-T_{i-1}(t), \; i=1,2,\ldots.
\end{equation}
Since these polynomials are orthogonal on the segment \([-1,1]\), we apply a linear mapping to the segment \([0,T]\).

The controlled norm of residual corresponding to the selected values of \(m,m_1\) and \(m_2\) takes the form
\begin{multline}\label{VSeries-8}
     \varepsilon_N=\max\limits_{t\in[0,T]}\left|y(t)-\sum\limits_{i=0}^{m-1}A_i\int\limits_{0}^{t}B_i(s)x(t-s)\,ds - \right. \\
     \left.- \sum\limits_{i=0}^{m_1-1}\sum\limits_{j=0}^{m_2-1}C_{ij}\int\limits_{0}^{t}\int\limits_{0}^{t}B_i(s_1)B_j(s_2)x(t-s_1)x(t-s_2)\,ds_1ds_2\right|
\end{multline}

Let us denote $N = m + m_1 m_2 - 1$. Number of equalities (number of nodes in the grid) equals to the number of unknown coefficients.

Thus, we have the following system of linear algebraic equations
\begin{equation}\label{VSeries-7}
  y_k=\sum\limits_{i=0}^{m-1}A_i\beta_{ik} + \sum\limits_{i=0}^{m_1-1}\sum\limits_{j=0}^{m_2-1}C_{ij}\gamma_{ijk}, \quad k=\overline{0,m + m_1 m_2 - 1},
\end{equation}
with respect to the unknown coefficients \(A_i,\;i=0,1,\ldots,m-1\) and \(C_{ij},\;i=0,1,\ldots,m_1-1,\;j=0,1,\ldots,m_2-1\). Here
\begin{equation}\label{VSeries-7.1}
  \beta_{ik}=\int\limits_{0}^{t_k}B_i(s)\,x(t_k-s)\,ds, \quad \gamma_{ijk}=\int\limits_{0}^{t_k}\int\limits_{0}^{t_k}B_i(s_1)B_j(s_2)x(t_k-s_1)x(t_k-s_2)\,ds_1ds_2.
\end{equation}

\section{Least--square method}
Let us denote $N > m + m_1m_2 - 1$. We have the situation where number of equalities is larger than number of unknown coefficients $A_i$ and $C_{ij}$. Thus we have the overdetermined system of linear equations with respect to the unknown coefficients \(A_i,\;i=0,1,\ldots,m-1\) and \(C_{ij},\;i=0,1,\ldots,m_1-1,\;j=0,1,\ldots,m_2-1\):
\begin{equation}
  y_k=\sum\limits_{i=0}^{m-1}A_i\beta_{ik} + \sum\limits_{i=0}^{m_1-1}\sum\limits_{j=0}^{m_2-1}C_{ij}\gamma_{ijk},
\end{equation}
where
\begin{equation}
  \beta_{ik}=\int\limits_{0}^{t_k}T_i(s)x(t_k-s)\,ds, \quad \gamma_{ijk}=\int\limits_{0}^{t_k}\int\limits_{0}^{t_k}T_i(s_1)T_j(s_2)x(t_k-s_1)x(t_k-s_2)\,ds_1ds_2.
\end{equation}

The system is inconsistent. Least--square method is used to find the approximate solution of the system. The point of the method is to find such coefficients $A_i$ and $C_{ij}$ such that the following criteria is minimized:
\begin{equation}
    \sum\limits_{k = 0}^{N - 1} \left( y_k - \sum\limits_{i=0}^{m-1}A_i\beta_{ik} - \sum\limits_{i=0}^{m_1-1}\sum\limits_{j=0}^{m_2-1}C_{ij}\gamma_{ijk} \right)^2 \longrightarrow \min
\end{equation}

\section{Numerical experiments}
Let us illustrate the operation of the proposed identification methods on two pairs of model signals.

\subsection{Model 1. Periodic signal}

\begin{equation}\label{VSeries-M1}
   \begin{split}
     x(t)=\sin(20t), \; y(t)=\frac{1}{81002}\Bigl(199\cos^2(20t)-15\sin(40t)-200\cos(20t)e^{-2t}+1+ \\
     +10\sin(20t)e^{-2t}+20\sin(20t)e^{-t}\Bigr) + \frac{1}{409}\left(3\sin(20t)-20\cos(20t)+\frac{850920}{40501}e^{-3t}\right).
   \end{split}
\end{equation}

\begin{figure}[h!]
   \center{\includegraphics[scale=0.6]{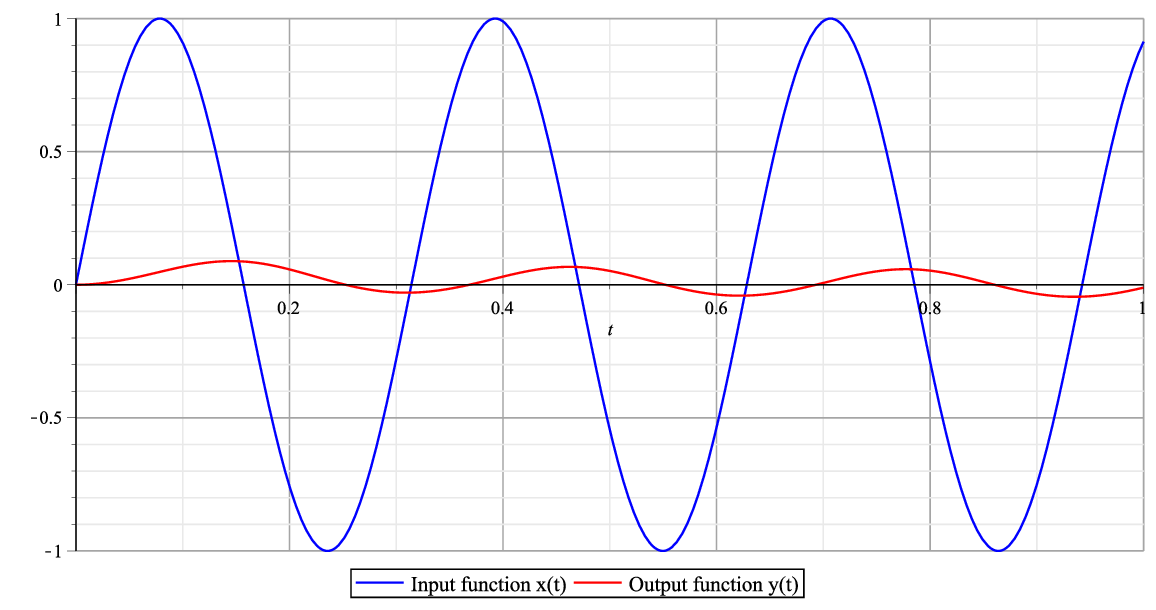}}
   \caption{Input and output functions \eqref{VSeries-M1}} \label{VSeries-P1}
\end{figure}

The Figure \ref{VSeries-P1} shows the graphs of the input and output signal \eqref{VSeries-M1}.

\subsubsection{Collocation method results for the model \eqref{VSeries-M1}}
The Table~\ref{VSeries-T1} demonstrates the dependence of the residual \(\varepsilon_N\) on the values \(m=m_1=m_2\) for the uniform mesh \(t_k=\frac{k}{N},\;k=0,1,\ldots,N,\) covering the segment \([0,1]\).
\begin{table}[H] 
\caption{Dependence of the residual \(\varepsilon_N\) on the values \(m,m_1,m_2\).}\label{VSeries-T1}
\newcolumntype{C}{>{\centering\arraybackslash}X}
\begin{tabularx}{\textwidth}{CC}
\toprule
$\mathbf{m}$	&  $\varepsilon_N$\\
\midrule
3		& $1,41\cdot 10^{-2}$			\\
4		& $1,14\cdot 10^{-6}$			\\
5		& $4,72\cdot 10^{-9}$		    \\
6		& $1,77\cdot 10^{-12}$			\\
7		& $1,83\cdot 10^{-14}$			\\
8		& $1,53\cdot 10^{-18}$			\\
10		& $2,84\cdot 10^{-26}$			\\
\bottomrule
\end{tabularx}
\end{table}

All calculations were performed in the Maple system with parameter \texttt{Digits:=30} (the number of digits that Maple uses when making calculations with software floating-point numbers).
Also note that the integration during the formation of the system \eqref{VSeries-7} was carried out analytically and did not introduce additional error in the calculation results. This is due to the fact that the input signal \(x(t)\) in most cases allows analytical calculation of the values \eqref{VSeries-7.1}. In the case of using input signals of a more complex structure, special approximation methods should be applied to the integrals \eqref{VSeries-7.1}, taking into account the possible fast oscillation of \(x(t)\).
\begin{figure}[H]
\vbox{
\begin{minipage}[h]{0.45\linewidth}
	\center{\includegraphics[width=1\linewidth]{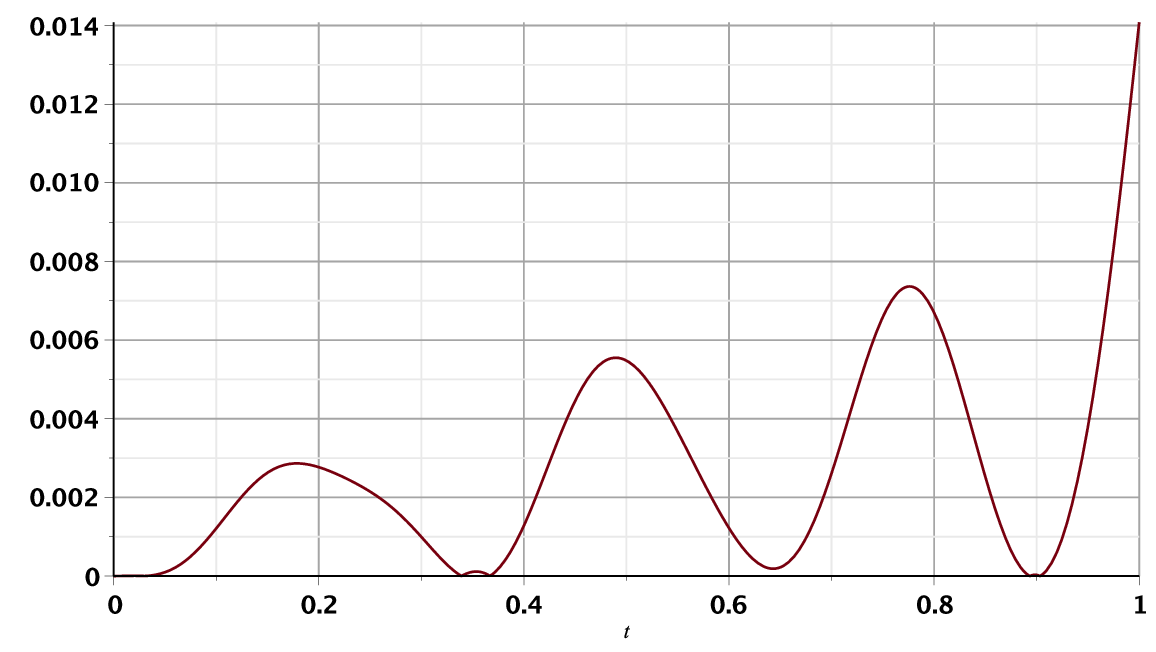}}
    \caption{Residual for \(m=3\)}
\end{minipage}
\begin{minipage}[!ht]{0.45\linewidth}
	\center{\includegraphics[width=1\linewidth]{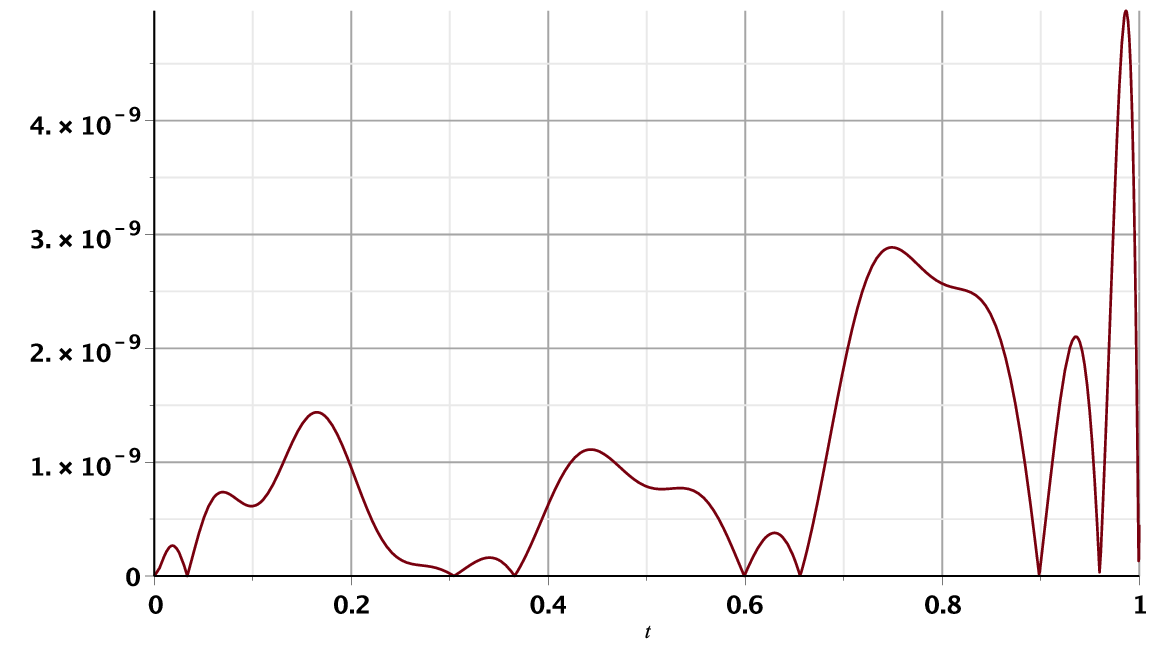}}
    \caption{Residual for \(m=5\)}
\end{minipage}\label{VSeries-P2}
}
\end{figure}
\begin{figure}[H]
\vbox{
\begin{minipage}[h]{0.45\linewidth}
	\center{\includegraphics[width=1\linewidth]{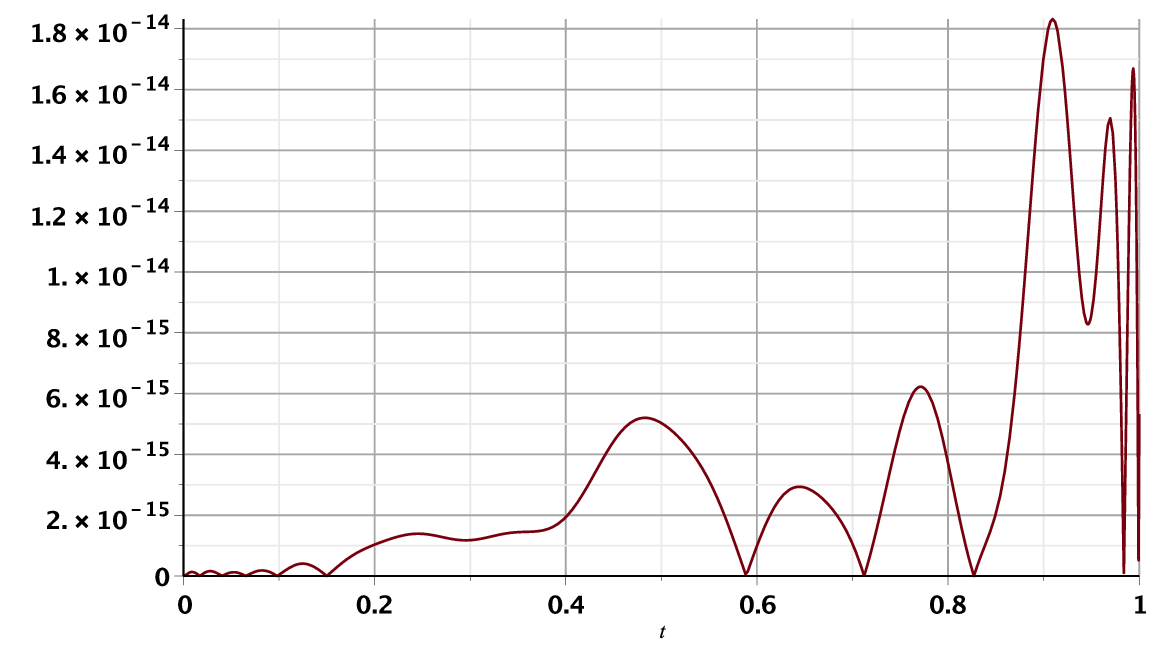}}
    \caption{Residual for \(m=7\)}
\end{minipage}
\begin{minipage}[!ht]{0.45\linewidth}
	\center{\includegraphics[width=1\linewidth]{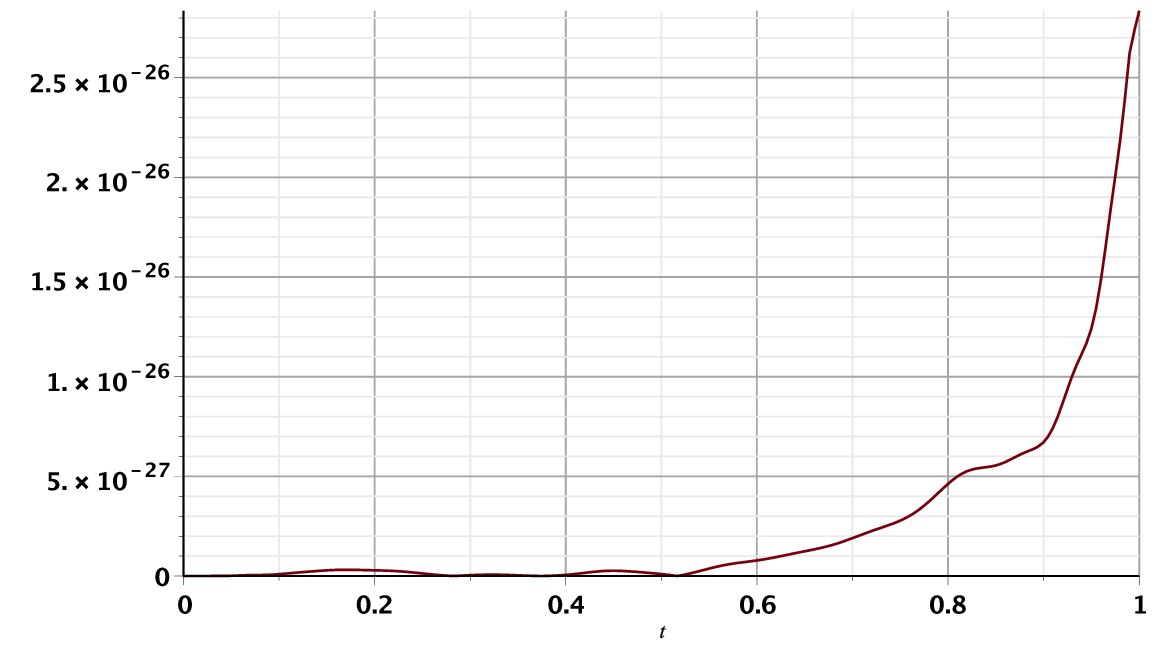}}
    \caption{Residual for \(m=10\)}
\end{minipage}\label{VSeries-P3}
}
\end{figure}

\subsubsection{Least--square method results for the model \eqref{VSeries-M1}}

For the simplicity we assume that $m = m_1 = m_2$. The Table~\ref{VSeries-T2} demonstrates the dependence of the residual \(\varepsilon_N\) on the parameters.

\begin{table}[H] 
\caption{Dependence of the residual \(\varepsilon_N\) on the values $m$ and $k$.\label{VSeries-T2}}
\newcolumntype{C}{>{\centering\arraybackslash}X}
\begin{tabularx}{\textwidth}{C|CCC}
\toprule
& \textbf{m = 3}	& \textbf{m = 5}	& \textbf{m = 7}\\
\midrule
$\mathbf{k = (m + m^2) \cdot 2}$ & $8.07 \cdot 10^{-4}$ & $4.92 \cdot 10^{-10}$ & $2.50 \cdot 10^{-16}$ \\
$\mathbf{k = (m + m^2) \cdot 5}$ &  $8.07 \cdot 10^{-4}$ & $3.90 \cdot 10^{-10}$ & $1.50 \cdot 10^{-16}$ \\ 
$\mathbf{k = (m + m^2) \cdot 10}$ & $8.07 \cdot 10^{-4}$ & $4.90 \cdot 10^{-10}$ &  $2.87 \cdot 10^{-15}$ \\
\bottomrule
\end{tabularx}
\end{table}

All calculations for least--square method were performed in MATLAB. Overdetermined matrix is solved using \texttt{lsqminnorm} function. It also should be noted that all the integrations during calculation were carried out analytically and didn't introduce additional error in the results.

\begin{figure}[H]
	\center{\includegraphics[width=1\linewidth]{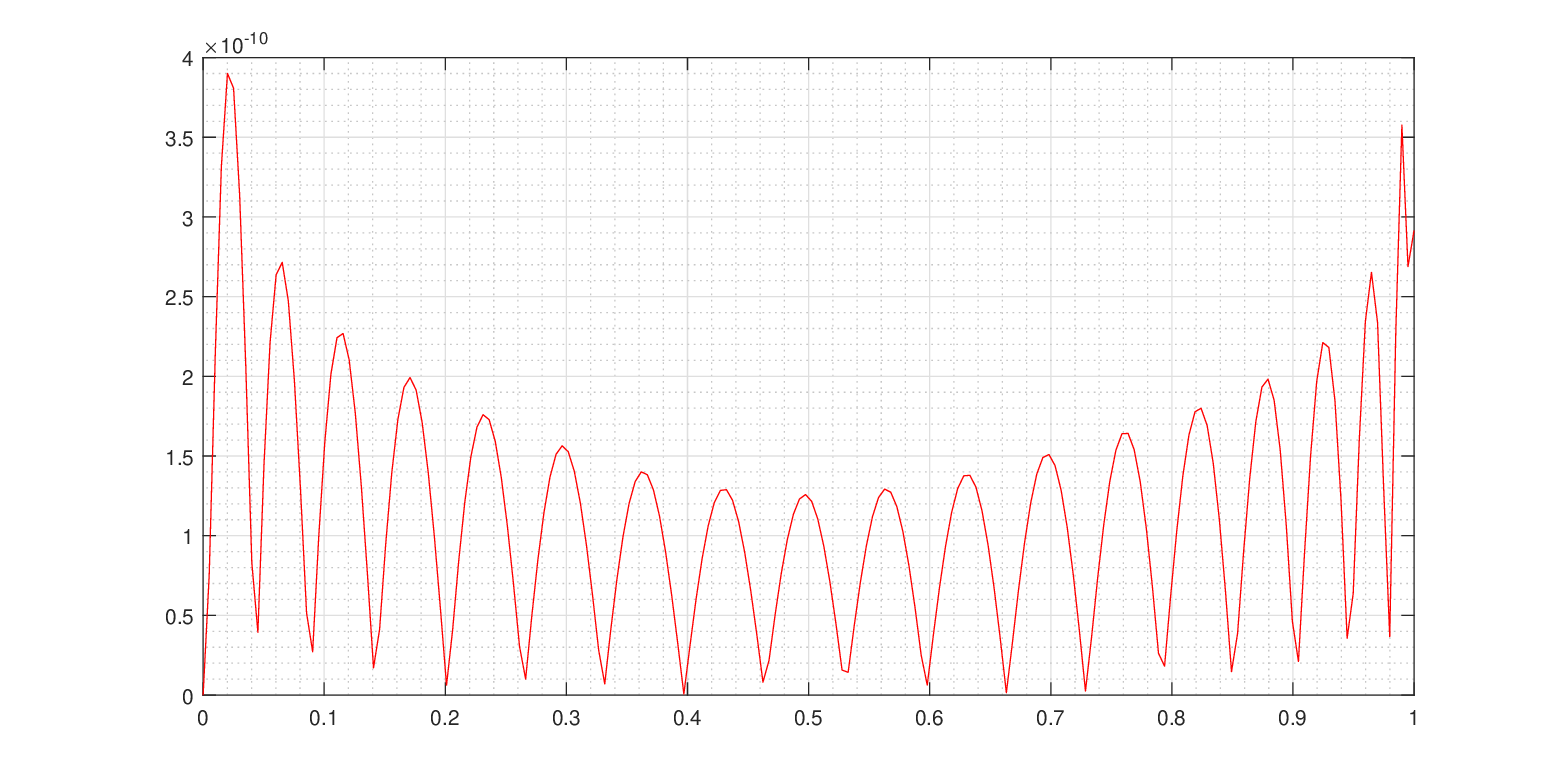}}
    \caption{Residual for $m=5$ and $k = (m + m^2) \cdot 5$}
\end{figure}
\begin{figure}[H]
	\center{\includegraphics[width=1\linewidth]{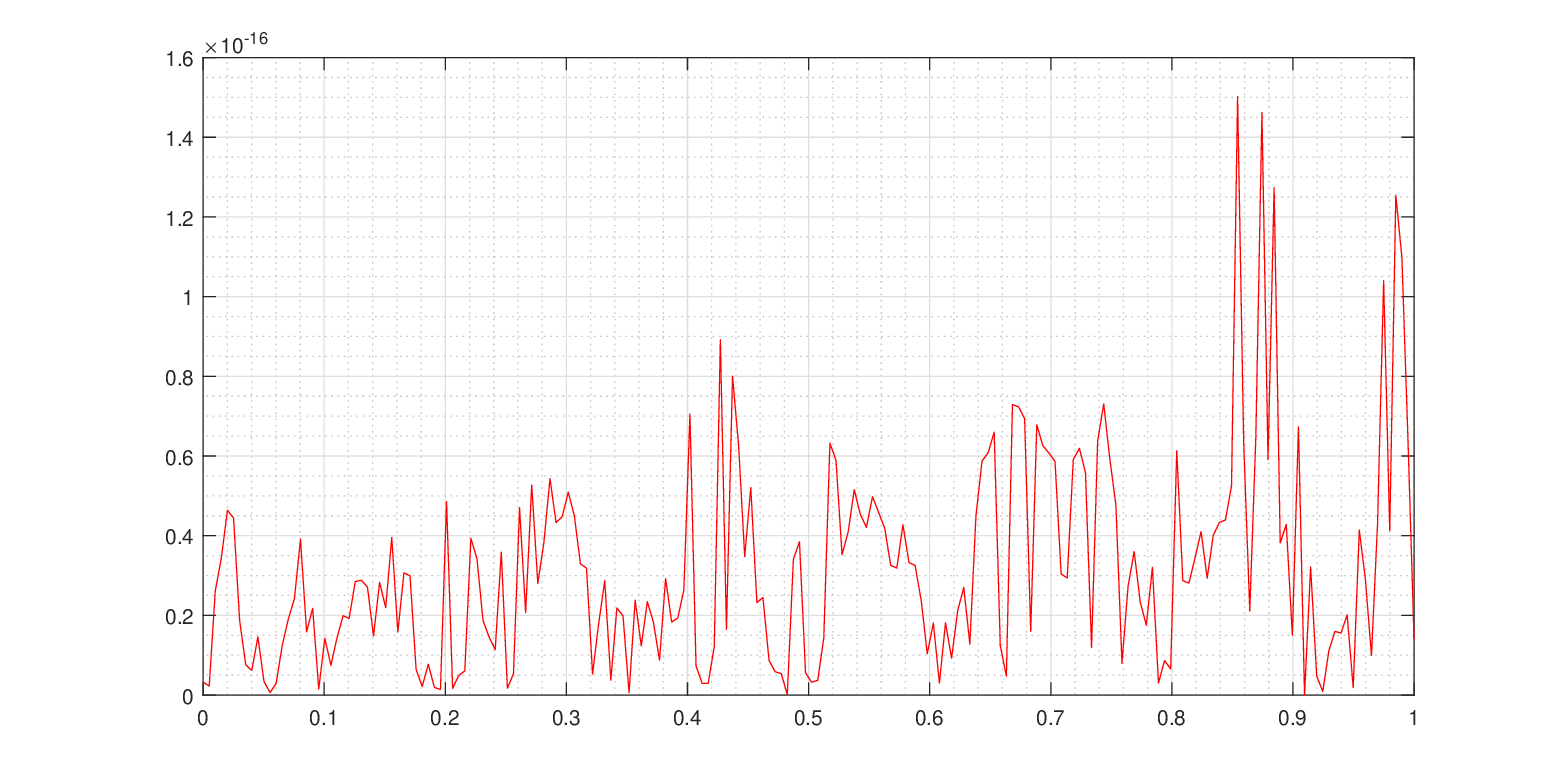}}
    \caption{Residual for $m=7$ and $k = (m + m^2) \cdot 5$}
\end{figure}

\subsection{Model 2. Fading input signal}
\begin{figure}[h!]
   \center{\includegraphics[scale=0.6]{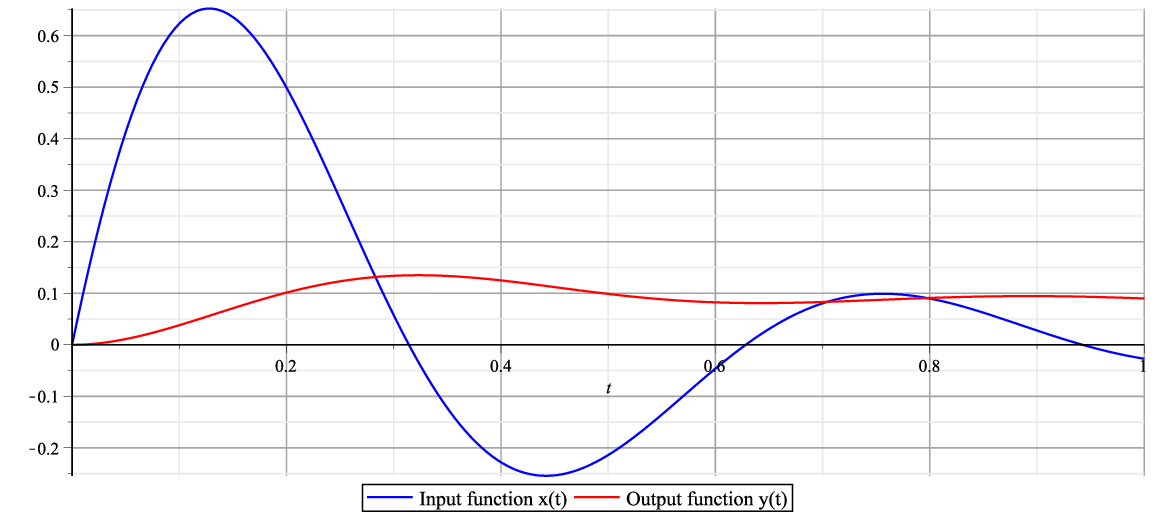}}
   \caption{Input and output functions \eqref{VSeries-M2}} \label{VSeries-P4}
\end{figure}

\begin{equation}\label{VSeries-M2}
   \begin{split}
     x(t)=e^{-3t}\sin(10t), \\ y(t)=\int\limits_{0}^{t}\cos\left(\frac{s}{2}\right)x(t-s)\,ds + \\ +\int\limits_{0}^{t}\int\limits_{0}^{t}\sin(s_1+2s_2)x(t-s_1)x(t-s_2)\,ds_1ds_2.
   \end{split}
\end{equation}
The Figure \ref{VSeries-P4} shows the graphs of the input and output signals \eqref{VSeries-M2}.

\subsubsection{Collocation method results for the model \eqref{VSeries-M2}}
The Table \ref{VSeries-T3} demonstrates the dependence of the residual \(\varepsilon_N\) on the values \(m=m_1=m_2\) for the uniform mesh \(t_k=\frac{k}{N},\;k=0,1,\ldots,N,\) covering the segment \([0,1]\).

\begin{table}[H] 
\caption{Dependence of the residual \(\varepsilon_N\) on the values \(m,m_1,m_2\).}\label{VSeries-T3}
\newcolumntype{C}{>{\centering\arraybackslash}X}
\begin{tabularx}{\textwidth}{CC}
\toprule
$\mathbf{m}$	&  $\varepsilon_N$\\
\midrule
3		& $3,16\cdot 10^{-5}$			\\
4		& $9,85\cdot 10^{-9}$			\\
5		& $8,58\cdot 10^{-12}$		    \\
6		& $2,17\cdot 10^{-16}$			\\
7		& $5,37\cdot 10^{-20}$			\\
\bottomrule
\end{tabularx}
\end{table}
\begin{figure}[!ht]
\vbox{
\begin{minipage}[h]{0.45\linewidth}
	\center{\includegraphics[width=1\linewidth]{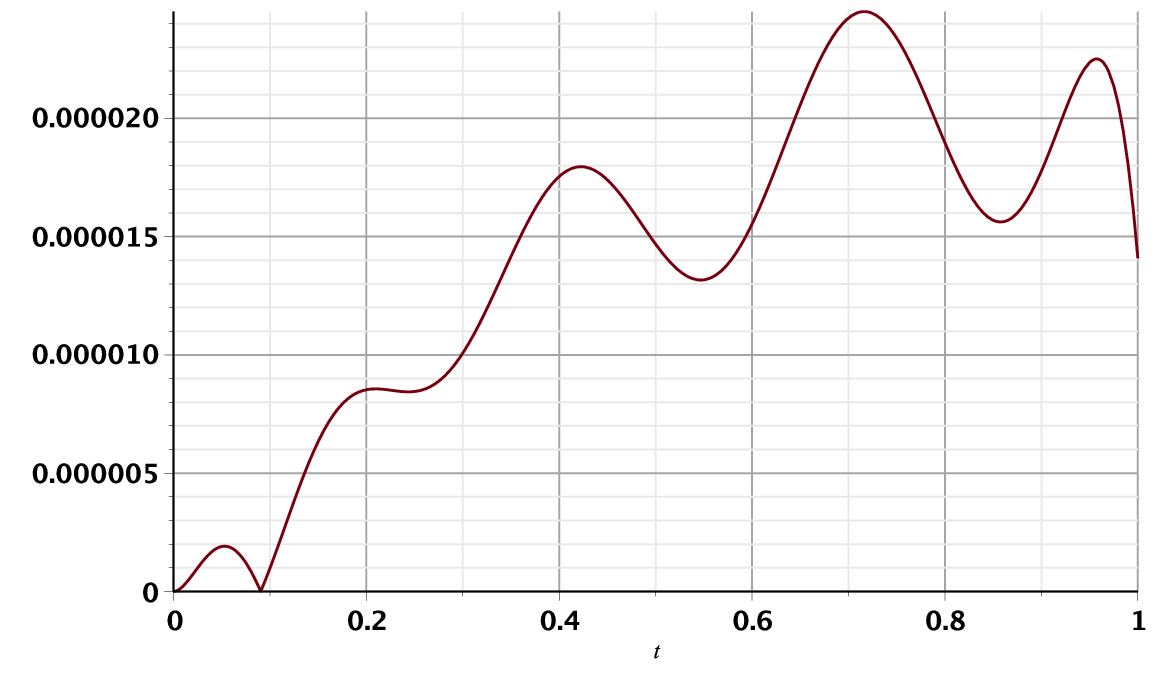}}
    \caption{Residual for \(m=3\)}
\end{minipage}
\begin{minipage}[!ht]{0.45\linewidth}
	\center{\includegraphics[width=1\linewidth]{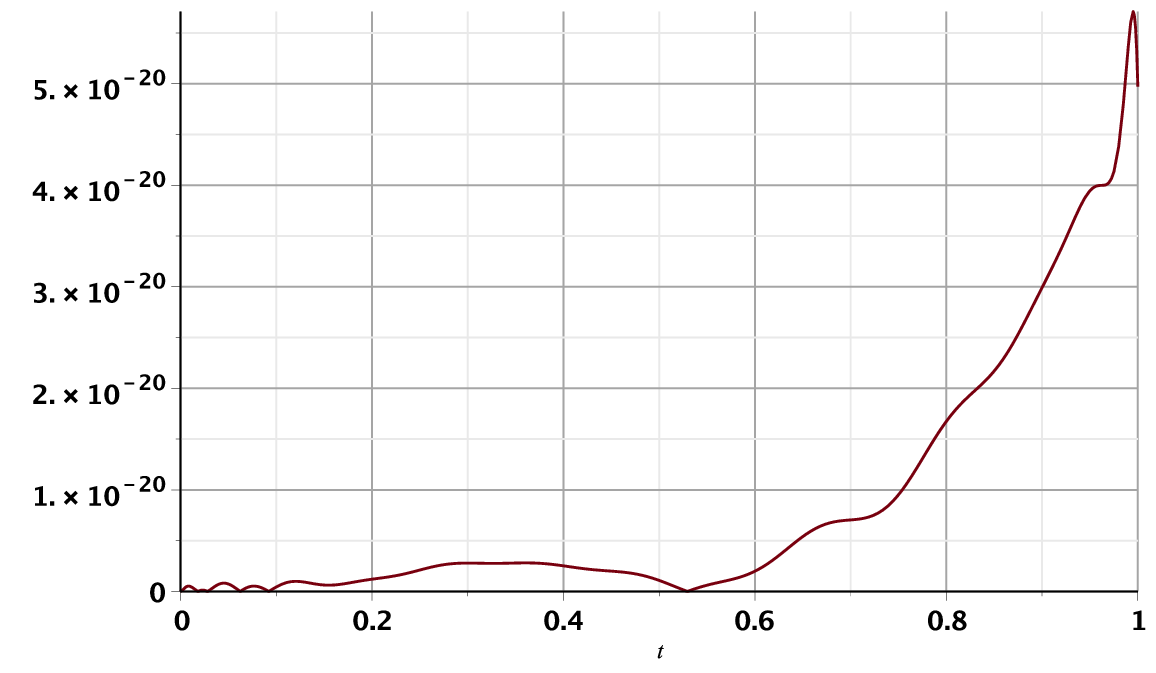}}
    \caption{Residual for \(m=7\)}
\end{minipage}\label{VSeries-P5}
}
\end{figure}

Let us also discuss the stability of suggested numerical technique. Let the input data of the problem \eqref{VSeries-M2} be determined with some random error \(\varepsilon_{rand}\) varying within the \(\delta\) value, namely \(\left|\varepsilon_{rand}\right|\leqslant\delta\).
Table \ref{VSeries-T4} shows the dependence of the averaged residual \(\varepsilon_N\) on the  \(\delta\) value at a fixed \(m=3\) based on the results of \(10\) measurements.
\begin{table}[H] 
\caption{Stability results for collocation}\label{VSeries-T4}
\newcolumntype{C}{>{\centering\arraybackslash}X}
\begin{tabularx}{\textwidth}{CC}
\toprule
$\delta$	&  $\varepsilon_N$\\
\midrule
$ 10^{-2}$		& $0.01729$			\\
$ 10^{-3}$		& $2.71\cdot 10^{-3}$			\\
$ 10^{-4}$		& $2.56\cdot 10^{-4}$		    \\
$ 10^{-5}$		& $7,54\cdot 10^{-5}$			\\
$ 10^{-6}$		& $1,66\cdot 10^{-5}$			\\
\bottomrule
\end{tabularx}
\end{table}
It can be seen from the results of the Table \ref{VSeries-T4} that residual continuously depends on the limits of random measurement errors of the input and output signals. Thus, we can conclude about the stability of the suggested method.

\subsubsection{Least--square method results for the model \eqref{VSeries-M2}}

\begin{table}[H] 
\caption{Dependence of the residual \(\varepsilon_N\) on the values $m$ and $k$.\label{VSeries-T5}}
\newcolumntype{C}{>{\centering\arraybackslash}X}
\begin{tabularx}{\textwidth}{C|CCC}
\toprule
& \textbf{m = 3}	& \textbf{m = 5}	& \textbf{m = 7}\\
\midrule
$\mathbf{k = (m + m^2) \cdot 2}$ & $2.38 \cdot 10^{-6}$ & $7.77 \cdot 10^{-14}$ & $2.93 \cdot 10^{-16}$ \\
$\mathbf{k = (m + m^2) \cdot 5}$ &  $2.63 \cdot 10^{-6}$ & $7.46 \cdot 10^{-14}$ & $3.05 \cdot 10^{-16}$ \\ 
$\mathbf{k = (m + m^2) \cdot 10}$ & $3.48 \cdot 10^{-6}$ & $7.41 \cdot 10^{-14}$ &  $3.80 \cdot 10^{-16}$ \\
\bottomrule
\end{tabularx}
\end{table}

\begin{figure}[H]
	\center{\includegraphics[width=1\linewidth]{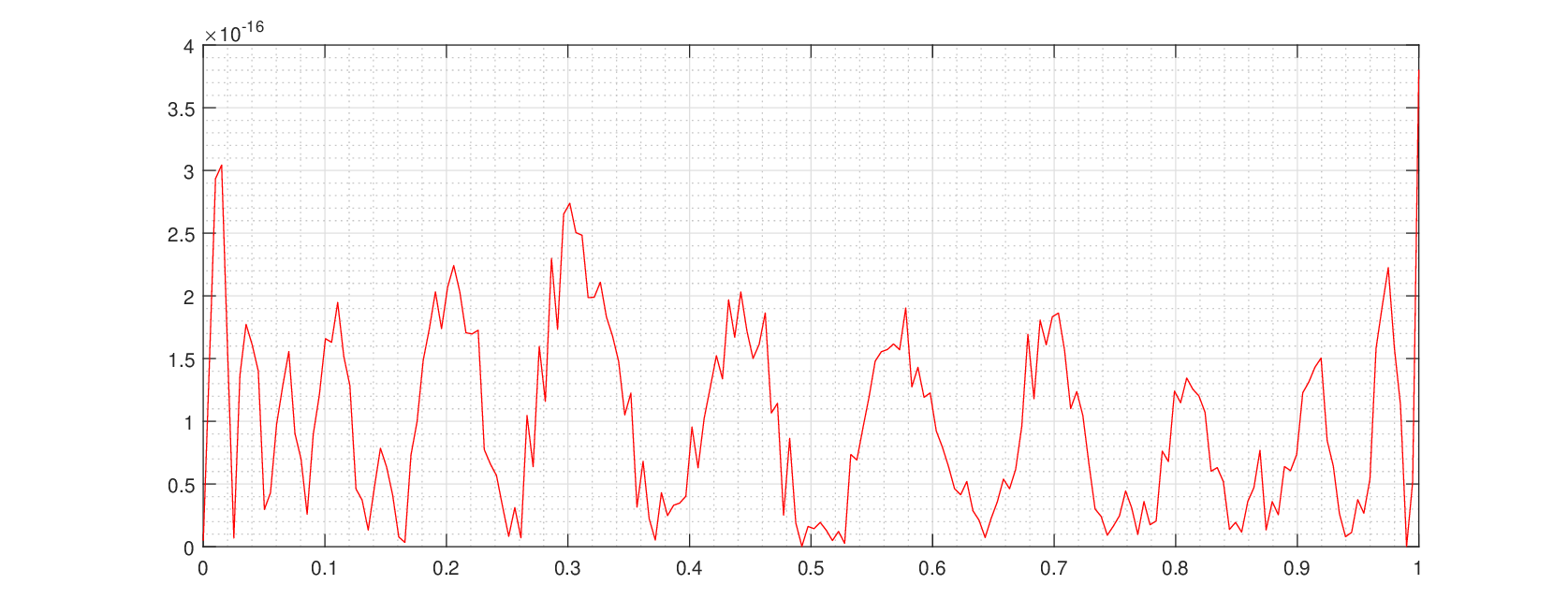}}
    \caption{Residual for $m=7$ and $k = (m + m^2) \cdot 7$}
\end{figure}

As for collocation method, let us check the stability of least-square method on this model. For testing stability 10 rounds of experiments were performed and average residual $\varepsilon_N$ was calculated. Also $m = 3, k = (m + m^2) \cdot 5$ were fixed. 

\begin{table}[H] 
\caption{Stability results for LSM}\label{VSeries-T6}
\newcolumntype{C}{>{\centering\arraybackslash}X}
\begin{tabularx}{\textwidth}{CC}
\toprule
$\delta$	&  $\varepsilon_N$\\
\midrule
$ 10^{-2}$		& $0.00628$			\\
$ 10^{-3}$		& $5.11 \cdot 10^{-4}$			\\
$ 10^{-4}$		& $6.02 \cdot 10^{-5}$		    \\
$ 10^{-5}$		& $5.36 \cdot 10^{-6}$			\\
$ 10^{-6}$		& $2.64 \cdot 10^{-6}$			\\
\bottomrule
\end{tabularx}
\end{table}


\section{Conclusions}
Two numerical approaches to solving the problem of identification of the Volterra model were proposed in the paper. As can be seen from the presented results, both methods showed stable convergence (in the sense of the tendency of the residual to zero). However, from the point of view of the arithmetic complexity of calculations, the collocation method turns out to be less expensive. And this factor is more pronounced the more parameters of the model are to be determined. This is due to the need to calculate a significantly larger number of integrals proportional to the square of the number of measurements being processed.

Further development of research suggests an increase in the number of terms \(n\) in the model \eqref{VS} to identify a more accurate functional relationship between the input and output signals. It is also planned to develop special methods for approximating integrals \eqref{VSeries-7.1} for the case of using input signals of a more complex structure, including fast oscillating signals.


\end{document}